# Elastic-plastic analysis of functionally graded bars under torsional loading

## George C. Tsiatas[a] and Nick G. Babouskos[b]

[a]*Corresponding author. e-mail:* gtsiatas@math.upatras.gr, *Department of Mathematics, University of Patras, Rio, GR-26504, Greece*

[b]*e-mail:* babouskosn@yahoo.gr, *School of Civil Engineering, National Technical University of Athens, Athens, GR-15773, Greece*

*Keywords*: Functionally graded materials, Non-homogeneous media, Elastic-plastic, Torsion.

**Abstract.** In this paper a new integral equation solution to the elastic-plastic problem of functionally graded bars under torsional loading is presented. The formulation is general in the sense that it can be applied to an arbitrary cross-section made of any type of elastoplastic material. In material science the Functionally Graded Material (FGM) is a non-homogeneous composite which performs as a single-phase material, by unifying the best properties of its constituent phase material. The nonlinear elastic-plastic behavior is treated by employing the deformation theory of plasticity. According to this theory, the material constants are assumed variable within the cross section, and are updated through an iterative process so as the equivalent stress and strain at each point coincide with the uniaxial material curve. In this investigation a new straightforward nonlinear procedure is introduced in the deformation theory of plasticity which simplifies the solution method. At each iteration step, the warping function is obtained by solving the torsion problem of a non-homogeneous isotropic bar using the Boundary Element Method (BEM) in conjunction with the Analog Equation Method (AEM). Without restricting the generality, the FGM material is comprised of a ceramic phase and a metal phase. The ceramic is assumed to behave linearly elastic, whereas the metal is modeled as an elastic – linear hardening material. Furthermore, the TTO homogenization scheme for estimating the effective properties of the two-phase FGM was adopted. Several bars with various cross-sections and material types are analyzed, in order to validate the proposed model and exemplify its salient features. Moreover, useful conclusion are drawn from the elastic-plastic behavior of functionally graded bars under torsional loading.

1. Introduction

In this paper a new integral equation solution to the elastic-plastic problem of functionally graded bars with arbitrary cross-section under torsional loading is presented. The continuous effort of engineers to design efficient materials which must be as light and economic as possible yet strong enough to withstand the most demanding functional requirements arising during their service life gave birth to a new class of materials; the Functionally Graded Material (FGM). In material science the FGM is a non-homogeneous composite which performs as a single-phase material, by unifying the best properties of its constituent phase materials. FGMs are deemed to have an advantageous behavior over laminated composites due to the continuous variation of their material properties yet in all three dimensions which alleviate delamination, de-bonding and matrix cracking initiation issues. There are several homogenization schemes for estimating the effective properties of a two-phase elastic FGM. Among them, the Mori-Tanaka [1] and self-consistent [2] schemes are the most prevalent ones. However, for nonlinear or elastic-plastic FGMs the problem is complicated due to additional material parameters. Besides the elastic moduli, which are quantified in linear FGMs, the flow stress and the plastic strain-hardening modulus must be also identified. Nakamura et al. [3] used the Kalman filter technique to estimate FGM through-thickness compositional variation and a rule-of-mixtures parameter that defines effective properties of FGMs. In this work the TTO homogenization scheme [4] is adopted for the evaluation of the FGM effective properties. Without restricting the generality, the FGM material is comprised of a ceramic phase and a metal phase. Bocciarelli et al. [5] extended the initially proposed TTO model to J2 flow theory with isotropic hardening to describe the elastoplastic behaviors of this metal-ceramic FGM. They concluded that the TTO model is an effective homogenization rule governing the transition from Hencky-Huber-Mises (HHM) model, typical of metals, toward a Drucker-Prager constitutive model which is more suitable to describe the mechanical response of ceramics. According to the TTO model the mechanical behavior of metal-ceramic composites beyond the elastic range, is essentially governed by the spreading of plasticity in the metal phase [5]. Specifically, the brittle ceramic is assumed to behave linearly elastic, whereas the metal is modeled as an elastic – linear hardening material.

Furthermore, a significant number of research papers has been published on the elastic-plastic torsion problem of homogeneous bars using analytical and numerical methods (see e.g. [6-20]). However, the work that has been conducted on the elastic-plastic torsion problem of composite bars is rather limited. Sapountzakis and Tsipiras [21] employed the Boundary Element Method (BEM) to the elastic–plastic uniform torsion problem of composite cylindrical bars of arbitrary cross-section consisting of materials in contact, each of which can surround a finite number of inclusions, taking into account the effect of geometric nonlinearity. Sapountzakis and Tsipiras [22] investigated also the effect of axial restraint on the previous elastic-plastic torsion problem of composite bars, treating the cases of free axial boundary conditions (vanishing axial force), restrained axial shortening or given axial force as special cases of an axially elastically supported bar. Recently, Bayat and Toussi [23] solved the elastoplastic torsion problem of hollow FGM circular shafts. The torsional shaft is considered as a thick-walled axisymmetric inhomogeneous cylindrical object, while the FG material is composed of ceramic and metallic parts with power function distribution only across the radial direction.

To the authors' knowledge, this is the first work that treats the elastic-plastic problem of functionally graded bars with arbitrary cross-section under torsional loading. The formulation is general in the sense that it can be applied to an arbitrary cross-section made of any type of elastoplastic material (elastic – perfectly plastic, strain hardening and nonlinear materials). The nonlinear elastic-plastic problem can be mathematically described by two general class-theories: (i) the total strain theory or deformation theory and (ii) the incremental strain theory or flow theory [24]. The former class-theory, known also as Hencky's theory, is employed in this investigation leading to the nonlinear elastoplastic boundary value problem. According to this theory the material constants are assumed variable within the cross section, and are updated through an iterative process so as the equivalent stress and strain at each point coincide with the uniaxial material curve. Several iterations schemes have been proposed in the literature which exhibit, however, a certain degree of complexity. In this investigation a new straightforward nonlinear procedure is introduced in the deformation theory of plasticity which simplifies the solution method. The key characteristic of the proposed approach lies in the fact that a nonlinear system of algebraic equations is constructed and any numerical method can be employed for its solution. At each iteration step, the warping function is obtained by solving the torsion problem of a non-homogeneous isotropic bar using the BEM in conjunction with the AEM [25], a robust integral equation method. Several bars with various cross-sections and material types are analyzed, in order to validate the proposed model and exemplify its salient features. Moreover, useful inferences are drawn from the elastic-plastic behavior of functionally graded bars under torsional loading.

## 2. Problem formulation

*2.1 Stress-strain relationship and effective material properties on total deformation theory*

Scrutinizing the total deformation theory of plasticity, it becomes evident that the formation of the constitutive relation between stress and strain is of vital importance. Following the work of Jahed et al. [26] this stress-strain relationship takes the form

$$\varepsilon_{ij} = f(\sigma_{ij}) \qquad (1)$$

where $f$ is a nonlinear function and $\varepsilon_{ij}$ is total strain tensor which is the sum of the conservative elastic and nonconservative plastic part [27]

$$\varepsilon_{ij} = \varepsilon_{ij}^e + \varepsilon_{ij}^p \qquad (2)$$

Using Hooke's law for isotropic material, the elastic strain tensor is related to the stress tensor as

$$\varepsilon_{ij}^e = \frac{1+\nu}{E}\sigma_{ij} - \frac{\nu}{E}\sigma_{kk}\delta_{ij} \qquad (3)$$

where $\nu$ is the Poisson's ratio, $E$ is the Young's modulus and $\delta_{ij}$ is the Kronecker delta.

According to Hencky's deformation theory, the plastic strain tensor is related to the deviatoric part of stress tensor in the following form [26-28]

$$\varepsilon_{ij}^p = \Phi s_{ij} \qquad (4)$$

where

$$s_{ij} = \sigma_{ij} - \frac{1}{3}\sigma_{kk}\delta_{ij} \qquad (5)$$

is the deviatoric stress tensor, and $\Phi$ is a scalar valued function given by

$$\Phi = \frac{3\varepsilon_{eq}^p}{2\sigma_{eq}} \qquad (6)$$

In the previous relation the equivalent plastic strain $\varepsilon_{eq}^p$ and equivalent stress $\sigma_{eq}$ are defined as

$$\varepsilon_{eq}^p = \sqrt{\frac{2}{3}\varepsilon_{ij}^p \varepsilon_{ij}^p} \qquad (7)$$

$$\sigma_{eq} = \sqrt{\frac{3}{2}s_{ij}s_{ij}} \qquad (8)$$

Substitution of Eqs. (4) - (8) into Eq. (3), provides

$$\varepsilon_{ij} = \left(\frac{1+\nu}{E} + \Phi\right)\sigma_{ij} - \left(\frac{\nu}{E} + \frac{\Phi}{3}\right)\sigma_{kk}\delta_{ij} \qquad (9)$$

The above equation can be rewritten as

$$\varepsilon_{ij} = \left(\frac{1+\nu_{eff}}{E_{eff}}\right)\sigma_{ij} - \frac{\nu_{eff}}{E_{eff}}\sigma_{kk}\delta_{ij} \qquad (10)$$

where $E_{eff}$ and $v_{eff}$ are the effective Young's modulus and Poisson's ratio, respectively, and are treated as material parameters. In general, Eq. (10) describes the material nonlinear behavior since $E_{eff}$ and $v_{eff}$ are field variables which uniquely depend on the final state of stress at every point of the material. Obviously, in the elastic regime the values of $E_{eff}$ and $v_{eff}$ are constants and equal to $E$ and $v$, respectively. Comparing Eqs. (9) and (10), we can obtain the effective values of the material parameters as

$$E_{eff} = \frac{1}{\frac{1}{E} + \frac{2}{3}\Phi} \tag{11}$$

$$v_{eff} = E_{eff}\left(\frac{v}{E} + \frac{\Phi}{3}\right) = \frac{1}{2} + \left(v - \frac{1}{2}\right)\frac{E_{eff}}{E} \tag{12}$$

Hencky suggested that nonlinear deformation can be modeled with strain-stress relations in which the final strains are a function of only the final stresses irrespective of the loading path to that final stress state [24]. In light of the above the deformation theory is applicable only if $E_{eff}$ and $v_{eff}$ can be expressed as a function of the multiaxial stress state. We have already expressed $v_{eff}$ in terms of $E_{eff}$ through Eq. (12). The final step is to express $E_{eff}$ as a function of the multiaxial stress state by use of the universal stress-strain curve.

The multiaxial stress state is described with the stress intensity [24]

$$\sigma_{eq} = \frac{\sqrt{2}}{2}\sqrt{(\sigma_x - \sigma_y)^2 + (\sigma_y - \sigma_z)^2 + (\sigma_z - \sigma_x)^2 + 6(\tau_{xy}^2 + \tau_{yz}^2 + \tau_{xz}^2)} \tag{13}$$

and the corresponding strain state is described with the strain intensity [24]

$$\varepsilon_{eq} = \frac{1}{\sqrt{2}(1+v_{eff})}\sqrt{(\varepsilon_x - \varepsilon_y)^2 + (\varepsilon_y - \varepsilon_z)^2 + (\varepsilon_z - \varepsilon_x)^2 + \frac{3}{2}(\gamma_{xy}^2 + \gamma_{yz}^2 + \gamma_{xz}^2)} \tag{14}$$

The above intensities are related to the effective Young's modulus as (see Fig. 1)

$$\sigma_{eq} = E_{eff}\varepsilon_{eq} \tag{15}$$

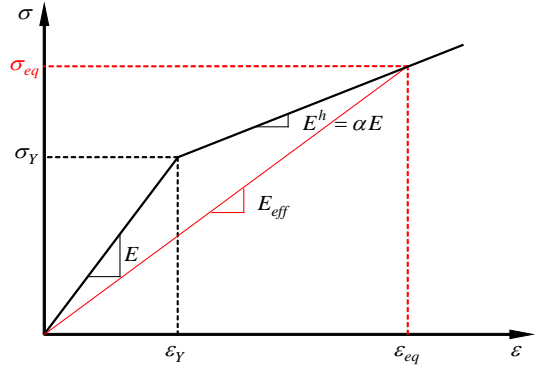

**Fig. 1.** The stress-strain relation for a linear work-hardening material.

Observing Fig. 1, it is deduced that $E_{eff}$ is the secant modulus which can be obtained by the uniaxial experimental material curve employing a nonlinear procedure described in the next section.

In the case of a linear work-hardening material (see Fig. 1) the stress-strain relationship can be cast in the form

$$\sigma = E\varepsilon \quad \text{for} \quad \varepsilon \leq \varepsilon_Y \tag{16}$$

$$\sigma = \sigma_Y(1-\alpha) + E^h \varepsilon \quad \text{for} \quad \varepsilon \geq \varepsilon_Y \tag{17}$$

where $\varepsilon_Y$ and $\sigma_Y$ are the yield strain and stress, respectively, and $E^h$ is the hardening modulus with $\alpha$ being the ratio of the slope of the linear hardening line to the slope of the elastic line. Thus, for the perfectly plastic case $\alpha = 0$ and, for the elastic case $\alpha = 1$.

*2.2 The torsion problem*

Consider now an elastic bar of length $L$ with arbitrary cross-section occupying the two-dimensional domain $\Omega$ of arbitrary shape in the $x, y$ plane bounded by the curve $\Gamma$ which may be piecewise smooth, i.e. it may have a finite number of corners. The cross-section is constant along the length of the bar and is twisted by moments $M_t$ applied at its ends. According to Saint-Venant's torsion theory (e.g. [14, 29]), the deformation of the bar consists of (a) rotations of the cross-sections about an axis passing through the *twist center* of the bar and (b) warping of the cross-sections, which is the same for all sections. Choosing the origin of the coordinate system at the twist center of an end section, the rotation at a distance $z$ is $\theta z$, where $\theta$ is a constant expressing the rotation of a cross-section per unit length. Assuming that this rotation is small, the displacement components of an arbitrary point are

$$u = -\theta z y \tag{18}$$

$$v = \theta z x \quad (19)$$

$$w = \theta \phi(x, y) \quad (20)$$

where $\phi(x, y)$ is the warping function. The displacement field given above yields the following nonzero components of the strain tensor

$$\gamma_{xz} = \theta(\phi_{,x} - y) \quad (21)$$

$$\gamma_{yz} = \theta(\phi_{,y} + x) \quad (22)$$

Substituting Eqs. (21) and (22) into the stress-strain relations of the deformation theory of plasticity Eqs. (10), yields the following nonzero components of the stress tensor

$$\tau_{xz} = G_{eff} \theta(\phi_{,x} - y) \quad (23)$$

$$\tau_{yz} = G_{eff} \theta(\phi_{,y} + x) \quad (24)$$

where

$$G_{eff} = \frac{E_{eff}}{2(1 + \nu_{eff})} \quad (25)$$

is the effective shear modulus. Recalling Eqs. (13) and (14) the stress $\sigma_{eq}$ and strain $\varepsilon_{eq}$ intensities for the Saint-Venant's torsion theory take the form

$$\sigma_{eq} = \sqrt{3(\tau_{yz}^2 + \tau_{xz}^2)} \quad (26)$$

$$\varepsilon_{eq} = \frac{1}{2(1 + \nu_{eff})} \sqrt{3(\gamma_{yz}^2 + \gamma_{xz}^2)} \quad (27)$$

Introducing Eqs. (23) and (24) into the equilibrium equations and the boundary conditions for the three-dimensional elastic body, noting that $\tau_{xz}$ and $\tau_{yz}$ are both independent of $z$, we obtain the following boundary value problem for the warping function [30]

$$G_{eff} \nabla^2 \phi + G_{eff,x} \phi_{,x} + G_{eff,y} \phi_{,y} = y G_{eff,x} - x G_{eff,y} \quad \text{in } \Omega \quad (28)$$

$$\nabla \phi \cdot \mathbf{n} = y n_x - x n_y \quad \text{on } \Gamma \quad (29)$$

where $\mathbf{n}$ is the unit vector, normal to the boundary $\Gamma$. Eq. (28) is actually the governing equation of the torsion problem of a non-homogeneous isotropic bar [30], since the effective shear modulus is treated as spatial field variable which varies at each point inside the domain $\Omega$.

Moreover, the torsional moment of the cross section is given as [29, 30]

$$M_t = \int_\Omega (x \tau_{yz} - y \tau_{xz}) d\Omega \quad (30)$$

## 2.3 FGM constitutive relation

In this investigation the TTO homogenization scheme [4] is adopted for the evaluation of the FGM effective properties. Without restricting the generality, the FGM material is comprised of a ceramic phase and a metal phase. In general, the volume fraction of the two FGM constituents (ceramic and metal) follows the power law distribution

$$V_c = (0.5 + y/h)^k, \quad V_c + V_m = 1 \tag{31, 32}$$

where $V$ is the volume fraction of the constituents, and $k$ is the non-negative power law exponent (see Fig. 2a). In our case, the multi-linear hardening elastoplastic material properties along the thickness can be defined by [3, 31]

$$E = (RE_m V_m + E_c V_c)/(RV_m + V_c), \quad R = \frac{q + E_c}{q + E_m} \tag{33, 34}$$

$$\nu = \nu_m V_m + \nu_c V_c, \quad \sigma_Y = \sigma_{Ym}\left(V_m + \frac{1}{R}\frac{E_c}{E_m}V_c\right) \tag{35, 36}$$

$$E^h = (RE_m^h V_m + E_c V_c)/(RV_m + V_c) \tag{37}$$

where $E$, $\nu$, $\sigma_Y$, and $E^h$ are the overall Young's modulus, Poisson's ratio, yield stress, and hardening modulus of the homogenized material (see Fig. 2b). Moreover, $q$ is the *stress transfer* dimensionless parameter [3]

$$q = \frac{(\sigma_c - \sigma_m)}{E_c(\varepsilon_c - \varepsilon_m)}, \quad (0 \leq q \leq \infty) \tag{38}$$

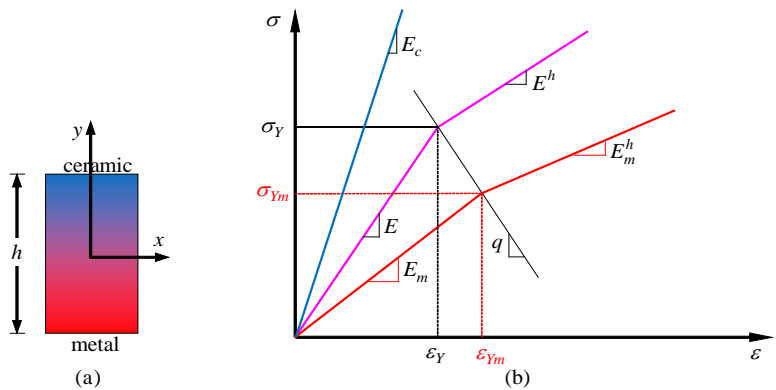

**Fig. 2.** Material distribution (a) and bilinear stress-strain relation for the FGM work-hardening material (b).

which expresses the normalized ratio of the stress to strain transfer. Generally, the parameter $q$ depends on several factors (e.g. mechanical characteristics of each constituent, loading etc.). However, in most applications $q$ is assumed constant even beyond the elastic limit due to lack of experimental data [5].

## 3. Numerical implementation

At each iteration step of the nonlinear procedure the warping function is obtained by solving the torsion problem of a non-homogeneous isotropic bar. In this section both the solution to the torsion problem and the nonlinear procedure are presented.

*3.1 The AEM solution to the non-homogeneous torsion problem*

The boundary value problem described by Eq. (28) together with the boundary condition Eq. (29), is solved using the BEM in conjunction with the AEM, a robust numerical integral equation method. According to the AEM, the original equation is substituted by an equation of the same order with known fundamental solution, which for the problem at hand is

$$\nabla^2 \phi = b(\mathbf{y}) \tag{39}$$

Eq. (39), which is called analog equation, indicate that the solution of Eq. (28) could be established by solving this Poisson's equation under the boundary condition Eq.(29), if the unknown source $b(\mathbf{y})$ ($\mathbf{y} \in \Omega$), was known. Its establishment is accomplished following the procedure below [30].

The integral representation of the solution to Eq. (39) is obtained by applying Green's identity for the warping function $\phi$ and $\phi^* = \ell n(r)/2\pi$ which is the fundamental solution to the Laplace equation. Thus, we have [29]

$$\varepsilon \phi(\mathbf{x}) = \int_\Omega \phi^*(\mathbf{x}, \mathbf{y}) b(\mathbf{y}) d\Omega_\mathbf{y} - \int_\Gamma [\phi^*(\mathbf{x}, \xi)\phi,_n(\xi) - \phi_i(\xi)\phi^*,_n(\mathbf{x}, \xi)] ds_\xi \tag{40}$$

in which $\mathbf{x} \in \Omega \cup \Gamma$, $\mathbf{y} \in \Omega$ and $\xi \in \Gamma$; $\phi^*,_n$ is the derivative of the fundamental solution normal to the boundary with $r = \|\mathbf{y} - \mathbf{x}\|$ or $r = \|\xi - \mathbf{x}\|$ being the distance between the points $\mathbf{x}, \mathbf{y}$ or $\mathbf{x}, \xi$. The expression (40) represents the solution of the differential equation (39) for points (i) inside the domain $\mathbf{x} \in \Omega$ ($\varepsilon = 1$), (ii) on the boundary $\mathbf{x} \in \Gamma$ ($\varepsilon = \alpha/2\pi$) where $\alpha$ is the angle between the tangents to the boundary at point $\mathbf{x}$ (for points where the boundary is smooth it is $\varepsilon = 1/2$), and (iii) outside the domain $\mathbf{x} \notin \Omega \cup \Gamma$ ($\varepsilon = 0$). For more information on the development of the BEM, we refer the interested reader to [29].

Applying now Eq. (40) to the boundary points yields the boundary integral equation

$$\frac{a}{2\pi}\phi(\mathbf{x}) = \int_{\Omega}\phi^*(\mathbf{x},\mathbf{y})b(\mathbf{y})d\Omega_{\mathbf{y}} - \int_{\Gamma}[\phi^*(\mathbf{x},\xi)\phi,_n(\xi) - \phi(\xi)\phi,_n^*(\mathbf{x},\xi)]ds_\xi \qquad (41)$$

which is a domain-boundary integral equation and could be solved using domain discretization to approximate the domain integrals. This, however, would spoil the advantages of BEM over the domain methods. We can maintain the pure boundary character of the method by converting the domain integrals to boundary line integrals using the following procedure. The unknown fictitious source $b(\mathbf{y})$ is approximated by the series

$$b(\mathbf{y}) = \sum_{j=1}^{M} a_j f_j \qquad (42)$$

where $f_j = f_j(r_j)$ is a set of $M$ radial basis approximation functions and $a_j$ are $M$ coefficients to be determined; $r_j = \|\mathbf{y} - \mathbf{x}_j\|$ with $\mathbf{x}_j$ being the $M$ collocation points located in the domain $\Omega$ (see Fig. 3). Using the Green's reciprocal identity [29] the domain integral in Eq. (41) becomes

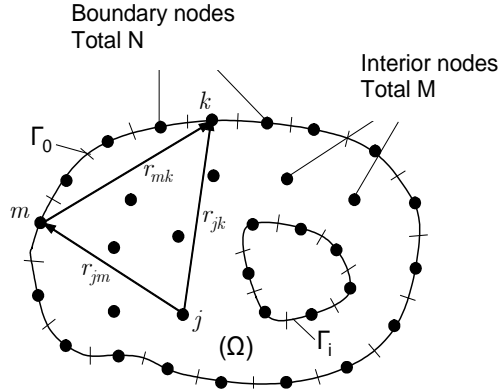

**Fig. 3.** Boundary discretization and domain nodal points.

$$\int_{\Omega}\phi^*(\mathbf{x},\mathbf{y})\phi(\mathbf{y})d\Omega_{\mathbf{y}} = \sum_{j=1}^{M} a_j \int_{\Omega}\phi^*(\mathbf{x},\mathbf{y})f_j(\mathbf{y})d\Omega_{\mathbf{y}}$$

$$= \sum_{j=1}^{M} a_j \left\{\varepsilon \hat{u}_j(\mathbf{x}) + \int_{\Gamma}\left[\phi^*(\mathbf{x},\xi)\hat{u}_j,_n(\xi) - \hat{u}_j(\xi)\phi,_n^*(\mathbf{x},\xi)\right]ds_\xi\right\} \qquad (43)$$

in which $\hat{u}_j(r_j)$ is a particular solution of the equation

$$\nabla^2 \hat{u}_j = f_j \qquad (44)$$

and $\hat{u}_j,_n$ its derivative normal to the boundary. Note that $\hat{u}_j$ can be always established when $f_j$ is specified. In view of Eq. (43), Eq. (41) is written as

$$\frac{a}{2\pi}\phi(\mathbf{x}) = \sum_{j=1}^{M} a_j \left\{ \varepsilon \hat{u}_j(\mathbf{x}) + \int_{\Gamma} \left[ \phi^*(\mathbf{x},\xi)\hat{u}_{j,n}(\xi) - \hat{u}_j(\xi)\phi^*_{,n}(\mathbf{x},\xi) \right] ds_\xi \right\}$$
$$- \int_{\Gamma} [\phi^*(\mathbf{x},\xi)\phi_{,n}(\xi) - \phi(\xi)\phi^*_{,n}(\mathbf{x},\xi)] ds_\xi \qquad (45)$$

For points $\mathbf{x} \in \Omega$ ($\varepsilon = 1$) the warping function $\phi(\mathbf{x})$ can be evaluated from Eq. (40), which by virtue of Eq. (43) provides

$$\phi(\mathbf{x}) = \sum_{j=1}^{M} a_j \left\{ \hat{u}_j(\mathbf{x}) + \int_{\Gamma} \left[ \phi^*(\mathbf{x},\xi)\hat{u}_{j,n}(\xi) - \hat{u}_j(\xi)\phi^*_{,n}(\mathbf{x},\xi) \right] ds_\xi \right\}$$
$$- \int_{\Gamma} [\phi^*(\mathbf{x},\xi)\phi_{,n}(\xi) - \phi(\xi)\phi^*_{,n}(\mathbf{x},\xi)] ds_\xi \qquad (46)$$

It is apparent that the first and second derivatives of the warping function $\phi(\mathbf{x})$ can be obtained by direct differentiation. Namely

$$\phi_{,kl}(\mathbf{x}) = \sum_{j=1}^{M} a_j \left\{ \hat{u}_{j,kl}(\mathbf{x}) + \int_{\Gamma} \left[ \phi^*_{,kl}(\mathbf{x},\xi)\hat{u}_{j,n}(\xi) - \hat{u}_j(\xi)\phi^*_{,nkl}(\mathbf{x},\xi) \right] ds_\xi \right\}$$
$$- \int_{\Gamma} [\phi^*_{,kl}(\mathbf{x},\xi)\phi_{,n}(\xi) - \phi(\xi)\phi^*_{,nkl}(\mathbf{x},\xi)] ds_\xi \qquad (47)$$

where $k, l = 0, x, y$. Note that $\phi_{,00} \equiv \phi$.

The final step of AEM is to apply Eq. (28) to the $M$ collocation points inside $\Omega$. We, thus, obtain a set of $M$ equations of the form

$$G^j_{eff} \nabla^2 \phi^j + G^j_{eff,x} \phi^j_{,x} + G^j_{eff,y} \phi^j_{,y} = y^j G^j_{eff,x} - x^j G^j_{eff,y} \qquad (48)$$

in which $j = 1, 2, ... M$. Substitution of Eqs. (47) into Eqs. (48) yield a set of $M$ linear algebraic equations which can be used to evaluate the coefficients $a_j$. This can be implemented only numerically using the procedure presented in the following.

The BEM with constant elements is used to approximate the boundary integrals in Eqs. (45), (46) and (47). In this case it is $a = \pi$, hence $\varepsilon = 1/2$. If $N$ is the number of the boundary nodal points (see Fig. 3), then for the $m$ nodal point Eq. (45) is written as

$$\frac{1}{2}\phi^m = \sum_{j=1}^{M} F_{mj} a_j + \sum_{k=1}^{N} \tilde{H}_{mk} \phi^k - \sum_{k=1}^{N} G_{mk} \phi^k_{,n} \qquad (49)$$

where

$$\tilde{H}_{mk} = \int_k \phi^*_{,n}(r_{mk}) ds, \quad G_{mk} = \int_k \phi^*(r_{mk}) ds \qquad (50), (51)$$

$$F_{mj} = \varepsilon \hat{u}^m_j - \sum_{k=1}^{N} \tilde{H}_{mk} \hat{u}^k_j + \sum_{k=1}^{N} G_{mk} (\hat{u}_{,n})^k_j \qquad (52)$$

in which $m = 1, 2, ..., N$ and $j = 1, 2, ..., M$.

Applying Eq. (49) to all boundary nodal points and using matrix notation provides

$$\mathbf{H}\boldsymbol{\phi} - \mathbf{G}\boldsymbol{\phi}_{,n} + \mathbf{F}\mathbf{a} = \mathbf{0} \tag{53}$$

where $\mathbf{a} = \{a_1, a_2, ..., a_M\}^T$; $\boldsymbol{\phi}$, $\boldsymbol{\phi}_{,n}$ are the vectors of the $N$ boundary nodal values of the warping function and its normal derivative, respectively, and

$$\mathbf{H} = \tilde{\mathbf{H}} - \tfrac{1}{2}\mathbf{I} \tag{54}$$

where $\mathbf{I}$ is the unit matrix.

The boundary condition Eq. (29) applied to all boundary points yields

$$\boldsymbol{\phi}_{,n} = \boldsymbol{\beta}_3 \tag{55}$$

where $\boldsymbol{\beta}_3$ is a known $N \times 1$ vector including the values of the right-hand-side of Eq. (29). Now, Eqs. (53) and (55) can be combined to express $\boldsymbol{\phi}$ in terms of $\mathbf{a}$

$$\boldsymbol{\phi} = \mathbf{S}\mathbf{a} + \mathbf{s} \tag{56}$$

in which $\mathbf{S}$ is known $N \times M$ rectangular matrix and $\mathbf{s}$ is known $N \times 1$ vector.

Moreover, using the same discretization in Eqs. (46), (47) and applying them to the $M$ collocation points we obtain

$$\tilde{\boldsymbol{\phi}} = \mathbf{F}\mathbf{a} + \mathbf{H}\boldsymbol{\phi} - \mathbf{G}\boldsymbol{\phi}_{,n} \tag{57}$$

$$\tilde{\boldsymbol{\phi}}_{,kl} = \mathbf{F}_{kl}\mathbf{a} + \mathbf{H}_{kl}\boldsymbol{\phi} - \mathbf{G}_{kl}\boldsymbol{\phi}_{,n} \tag{58}$$

where $k, l = 0, x, y$. In the previous expression $\mathbf{H}$, $\mathbf{H}_{kl}$, $\mathbf{G}$, $\mathbf{G}_{kl}$ are known $M \times N$ matrices originating from the integration of the kernels on the boundary elements of Eqs. (46), (47) and $\tilde{\boldsymbol{\phi}}$, $\tilde{\boldsymbol{\phi}}_{,kl}$ are $M \times 1$ vectors including the values of $\phi$ and its derivatives at the interior collocation points. Note that the resulting line integrals are regular since the distance $r_{jm}$ in the kernels does not vanish.

Subsequently, substituting Eqs. (58) into Eqs. (48) and using Eqs. (55) and (56) to eliminate $\phi$ and $\phi_{,n}$, yields the following linear system of equations for $\mathbf{a}$

$$\mathbf{K}\mathbf{a} = \mathbf{g} \tag{59}$$

Evidently, the coefficients $\mathbf{a}$ are employed in Eq. (53) to evaluate $\boldsymbol{\phi}$ while the warping function and its derivatives at the $M$ nodal points are computed from Eqs. (57) and (58). For points $P \in \Omega$ not coinciding with the nodal points these quantities are evaluated from the discretized counterparts of Eqs. (46) and (47).

The shear stresses are now evaluated using Eq. (23) and (24). The torsional moment of the cross section is obtained by the domain integral of Eq. (30). In this work we evaluate the domain integral using the method presented in [32], which transforms the domain integral into a boundary line one. Thus, the integrand $R(\mathbf{x}) = x\tau_{yz} - y\tau_{xz}$ of the domain integral is approximated by radial basis functions

$$R(\mathbf{x}) = \sum_{j=1}^{M} \bar{a}_j f_j \quad (60)$$

where $\bar{a}_j$ are $M$ coefficients to be determined. Collocating the above equation at the $M$ domain nodal points we obtain

$$\mathbf{R} = \mathbf{\Phi}\bar{\mathbf{a}} \quad (61)$$

where $\mathbf{\Phi} = [f_{ji}]$ is an $M \times M$ matrix and $\mathbf{R}$ is a $M \times 1$ vector with the values of the integrand at the domain nodal points. The $\bar{a}_j$ coefficients are obtained from Eq. (61) as

$$\bar{\mathbf{a}} = \mathbf{\Phi}^{-1}\mathbf{R} \quad (62)$$

Substituting Eq. (60) in Eq.(30) yields

$$M_t = \int_\Omega R(\mathbf{x})d\Omega = \int_\Omega \sum_{j=1}^{M} \bar{a}_j f_j d\Omega = \sum_{j=1}^{M} \bar{a}_j \int_\Omega f_j d\Omega \quad (63)$$

Using Green's reciprocal identity and the same boundary discretization as previous, Eq. (63) becomes

$$M_t = \sum_{j=1}^{M}\sum_{k=1}^{N} \bar{a}_j \int_{\Gamma_k} \hat{u}_{j,n}\, ds \quad (64)$$

which permits the evaluation of the torsional moment of the cross section.

*3.2 Nonlinear procedure*

The nonlinear procedure consists in finding the complete spatial distribution of $E_{\mathit{eff}}$ and $\nu_{\mathit{eff}}$ so that the equivalent stress and strain at each point must coincide with the uniaxial material curve. Several iterations schemes have been proposed in the literature for the determination of $E_{\mathit{eff}}$ and $\nu_{\mathit{eff}}$ (e.g. the projection method, the arc-length method, the Neuber's rule) which exhibit, however, a certain degree of complexity.

In this investigation a new straightforward nonlinear procedure is introduced in the deformation theory of plasticity which simplifies the solution method. The key characteristic of the proposed approach lies in the fact that a nonlinear system of algebraic equations is constructed and any numerical method can be employed for its solution.

To begin with, $E_{\mathit{eff}}$ is approximated by the series

$$E_{eff}(\mathbf{y}) = \sum_{j=1}^{M} k_j f_j \qquad (65)$$

where $f_j = f_j(r_j)$ is a set of $M$ radial basis approximation functions and $k_j$ are $M$ coefficients to be determined; $r_j = \|\mathbf{y} - \mathbf{x}_j\|$ with $\mathbf{x}_j$ being the $M$ collocation points located in the domain $\Omega$ (see Fig. 3).

Secondly, $v_{eff}$ and $G_{eff}$ are calculated by the following relations

$$v_{eff} = \frac{1}{2} + \left(v - \frac{1}{2}\right)\frac{E_{eff}}{E}, \quad G_{eff} = \frac{E_{eff}}{2(1+v_{eff})} \qquad (66)$$

Further, for a given value of $\theta$, the shear strains

$$\gamma_{xz} = \theta(\phi_{,x} - y), \quad \gamma_{yz} = \theta(\phi_{,y} + x) \qquad (67), (68)$$

and the accompanying stresses

$$\tau_{xz} = G_{eff}\theta(\phi_{,x} - y), \quad \tau_{yz} = G_{eff}\theta(\phi_{,y} + x) \qquad (69), (70)$$

are employed for the evaluation of the equivalents stress $\sigma_{eq}$ and strain $\varepsilon_{eq}$ intensities

$$\sigma_{eq} = \sqrt{3(\tau_{yz}^2 + \tau_{xz}^2)}, \quad \varepsilon_{eq} = \frac{1}{2(1+v_{eff})}\sqrt{3(\gamma_{yz}^2 + \gamma_{xz}^2)} \qquad (71), (72)$$

Consequently, Eqs. (69) and (70) are used to cast the system of the $M$ nonlinear algebraic equations for the unknown parameters $k_j$

$$\sigma_{eq}(k_j) = \begin{cases} E\varepsilon_{eq}(k_j) & , \text{ if } \varepsilon_{eq}(k_j) < \varepsilon_Y \\ \sigma_Y + E^h(\varepsilon_{eq}(k_j) - \varepsilon_Y), & \text{ if } \varepsilon_{eq}(k_j) < \varepsilon_Y \end{cases} \qquad (73)$$

in order to satisfy the universal material curve at each of the $M$ collocation points located in the domain $\Omega$.

Note that in Eqs. (67) - (70) the warping function and its derivatives have been evaluated by the solution of the torsion problem of a non-homogeneous isotropic bar described in the previous section.

The nonlinear procedure presented in this section permits the analysis of any type of elastoplastic material (elastic – perfectly plastic, strain hardening and nonlinear materials). In this work the nonlinear system of algebraic equations was solved using the Newton-Raphson iteration method.

## 4. Numerical examples

On the basis of the numerical procedure presented in the previous section, a FORTRAN code has been written and numerical results for bars with various cross-sections and material types have been obtained, which illustrate the applicability, effectiveness and accuracy of the proposed method. The employed radial basis functions $f_j$ are the multiquadrics (MQs) which are defined as $f_j = \sqrt{r^2 + c^2}$, where $c$ is a shape parameter. The particular solution $\hat{u}_j(r_j)$ is

$$\hat{u}_j = -\frac{c^3}{3}\ln\left(c\sqrt{r^2+c^2}+c^2\right)+\frac{1}{9}\left(r^2+4c^2\right)\sqrt{r^2+c^2} \tag{74}$$

In all cases the shape parameter was taken $c = 0.1$.

### 4.1 Elastic-perfectly plastic bar with rectangular cross-section

As a first example we treat the torsion problem of an isotropic elastic-perfectly plastic bar with rectangular cross-section $b \times h$ with $b = 5cm$ and $h = 10cm$. The employed data of the material are: $E = 210600 \text{ kN/cm}^2$, $v = 0.3$ and $\sigma_Y = 24 \text{ kN/cm}^2$. As the results depend mostly on the number of domain nodal points, a convergence test was carried out. Table 1 presents the torsional moment for various values of the domain nodal points $M$ using $N = 300$ boundary elements. It can be easily seen that 300 nodal points are sufficient in order to gain convergence solution even for large values of the rotation $\theta$. Nonetheless, in this study 450 domain nodal points were used in order to depict accurately the plastic regions of the cross section (see Fig 4). Additionally, the proposed method proves its computational efficiency since a few seconds are needed on a personal computer with an Intel Core i3 processor to get accurate results.

**Table 1**

Torsional moment for various values of $M$ ($N = 300$).

| $\theta/\theta_{el}$ | $M_t/M_{el}$ | | | | |
|---|---|---|---|---|---|
| | $M = 98$ | $M = 162$ | $M = 200$ | $M = 300$ | $M = 450$ |
| 1.09 | 1.04 | 1.06 | 1.07 | 1.07 | 1.08 |
| 1.50 | 1.38 | 1.37 | 1.37 | 1.36 | 1.36 |
| 1.90 | 1.53 | 1.52 | 1.51 | 1.51 | 1.50 |
| 2.45 | 1.64 | 1.62 | 1.62 | 1.60 | 1.58 |
| 3.00 | 1.66 | 1.65 | 1.65 | 1.64 | 1.63 |

The maximum value of the elastic torsional moment for this cross-section is $M_{el} = 0.142hb^2\sigma_Y = 852.2$ kNcm and the maximum plastic torsional moment is $M_{pl} = 0.0962(3h-b)b^2\sigma_Y = 1443.4$ kNcm [14]. Fig. 5 presents the torsional moment-rotation curve obtained by the presented method and the results are in very good agreement with that obtained by [21] using flow theory and the BEM. Fig. 6 shows the plastic zones of the cross section for various values of the rotation $\theta$. Moreover, in Fig. 7 the shear stress vector distribution in the elastic state ($\theta/\theta_{el}=1$) and in elastic-plastic state ($\theta/\theta_{el}=4.75$) are depicted.

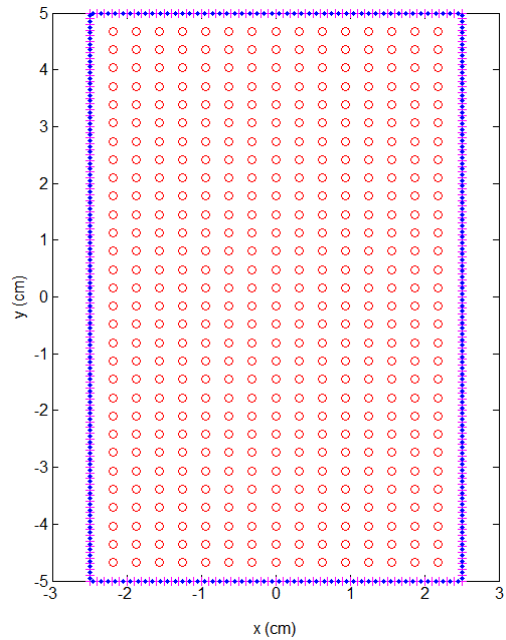

**Fig. 4.** Distribution of $N=300$ boundary and $M=450$ domain nodal points in Example 4.1.

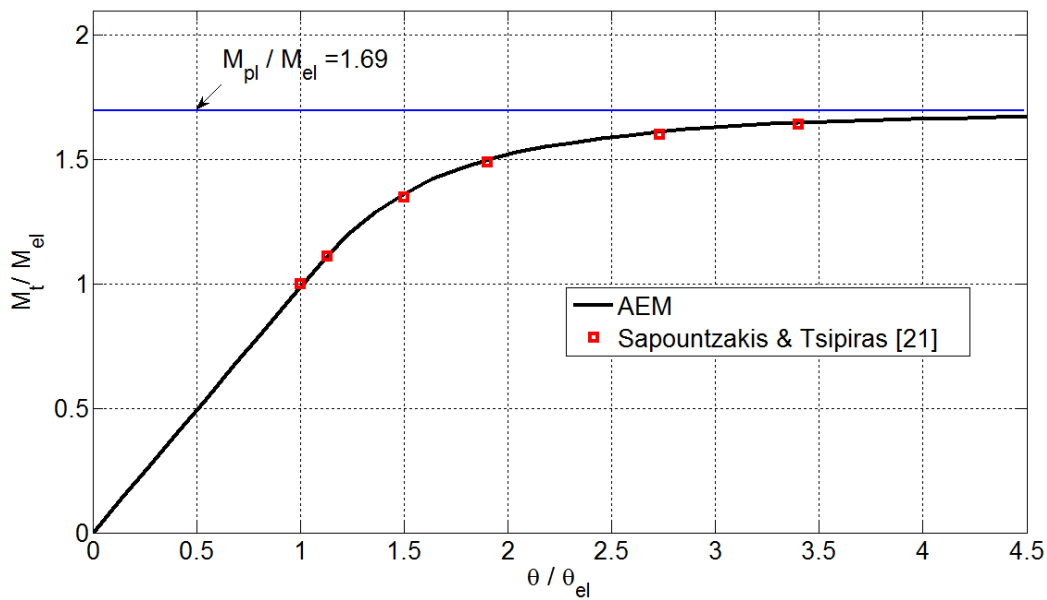

**Fig. 5.** Torsional moment-rotation curve in Example 4.1.

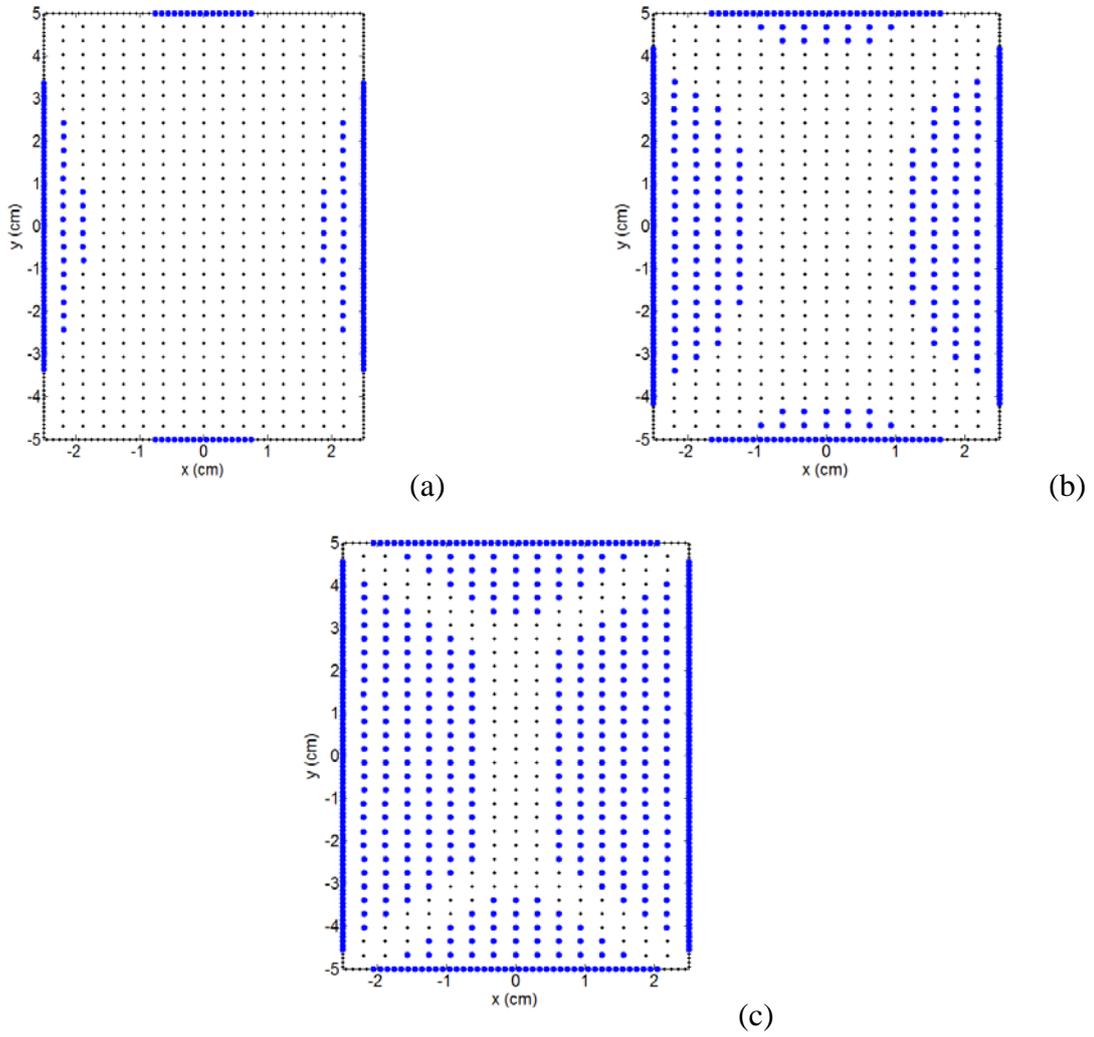

**Fig. 6.** Plastic regions of the rectangular cross section for (a) $\theta/\theta_{el}=1.36$ (b) $\theta/\theta_{el}=2.18$ and (c) $\theta/\theta_{el}=4.75$ in Example 4.1.

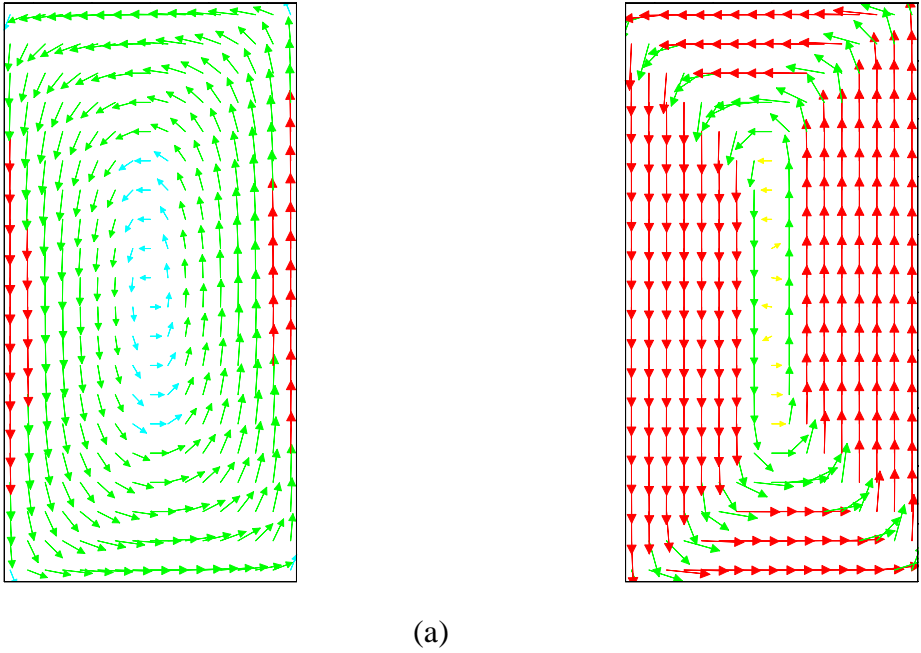

(a)                          (b)

**Fig. 7.** Shear stress ($kN/cm^2$) vector distribution in (a) elastic state ($\theta/\theta_{el}=1$) and (b) in elastic-plastic state at $\theta/\theta_{el}=4.75$ in Example 4.1.

*4.2 Elastic-plastic bar with triangular cross-section*

As a second example we study the torsion problem of an isotropic elastic-plastic bar with an equilateral triangular cross section with side length $b=10$ cm, as shown in Fig. 8. The employed data of the material are: $E=210600$ $kN/cm^2$, $\nu=0.3$ and $\sigma_Y=24$ $kN/cm^2$. The results were obtained using $N=240$ constant boundary elements and $M=288$ nodal points (see Fig. 8). The exact values of the maximum elastic and plastic torsional moments are given as: $M_{el}=0.02884\sigma_Y b^3=692.3$ kNcm and $M_{pl}=0.04811\sigma_Y b^3=1154.7$ kNcm. Fig. 9 shows the torsional moment-rotation curve for elastic-perfectly plastic material, which has been obtained by the presented method. The results are in very good agreement with that obtained by flow theory using the FEM [14]. More specifically, at $\theta/\theta_{el}=4$ the FEM in [14] gives $M_t/M_{el}=1.622$, whereas our proposed method gives $M_t/M_{el}=1.645$. Furthermore, Fig. 10 shows the plastic zones of the triangular cross-section for various values of the rotation $\theta$. Fig. 11 and 12 present the contours of the warping surface and the shear stress vector distribution in the elastic ($\theta/\theta_{el}=1$) and elastic-plastic regime ($\theta/\theta_{el}=4.3$), respectively. Finally, in Fig 13 the moment-rotation curves for hardening material with $E_h=aE$, for two cases of the parameter $a$ (0.3 and 0.5) are depicted.

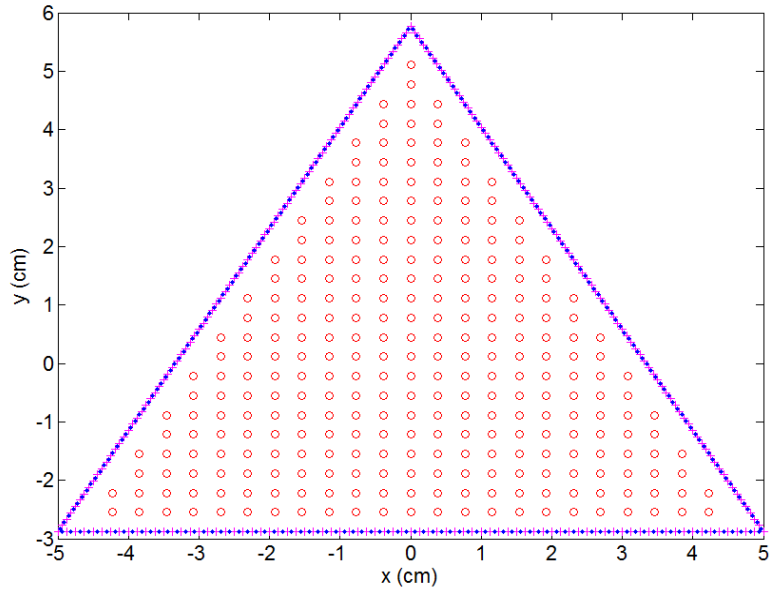

**Fig. 8.** Boundary and domain nodal points in triangular cross section in Example 4.2.

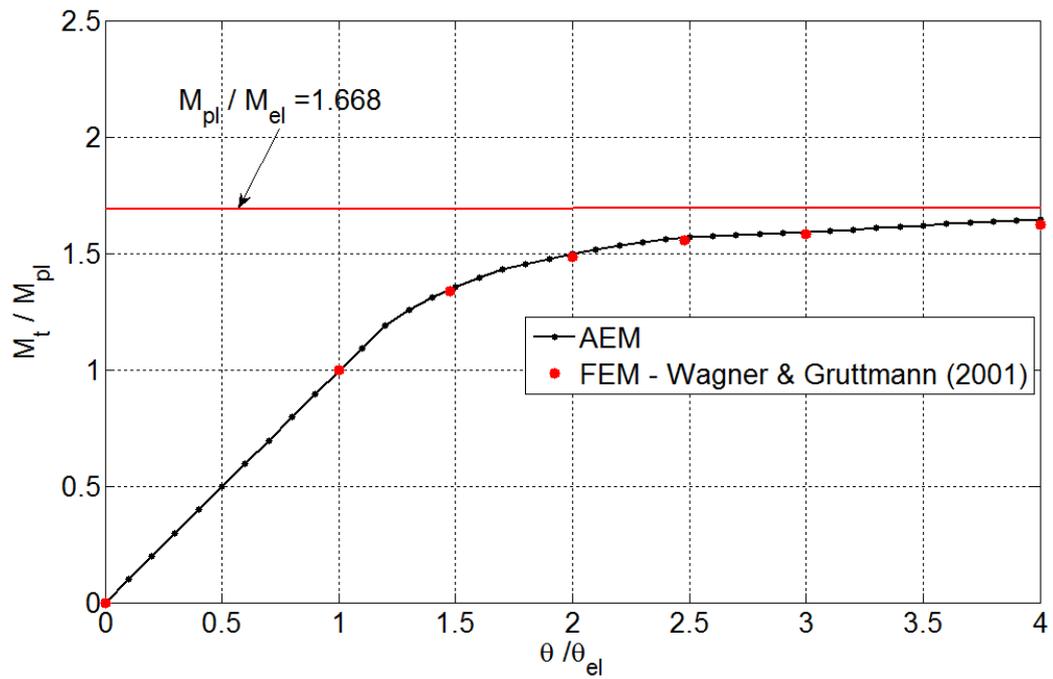

**Fig. 9.** Torsional moment-rotation curve in Example 4.2.

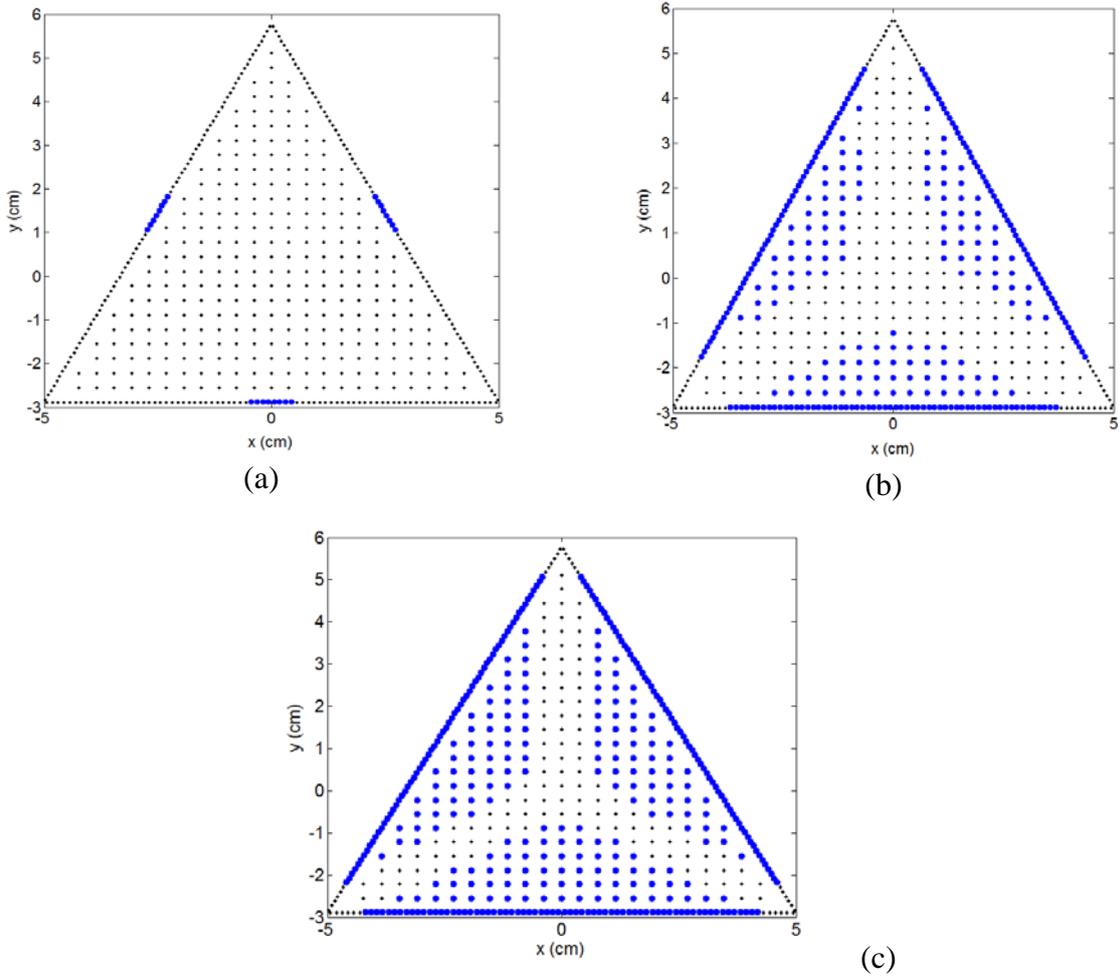

**Fig. 10.** Plastic regions of the triangular cross section for (a) $\theta/\theta_{el} = 1.05$ (b) $\theta/\theta_{el} = 2.6$ and (c) $\theta/\theta_{el} = 4.3$ in Example 4.2.

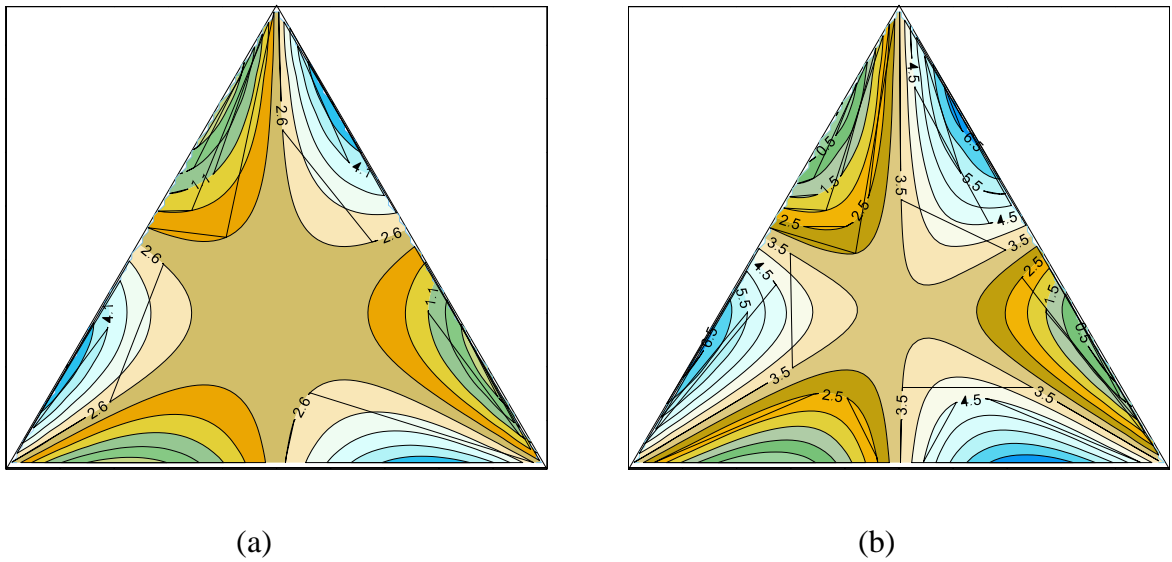

(a) (b)

**Fig. 11.** Contours of the warping surface in the (a) elastic ($\theta/\theta_{el} = 1$) and (b) elastic-plastic ($\theta/\theta_{el} = 4.3$) regime in Example 4.2.

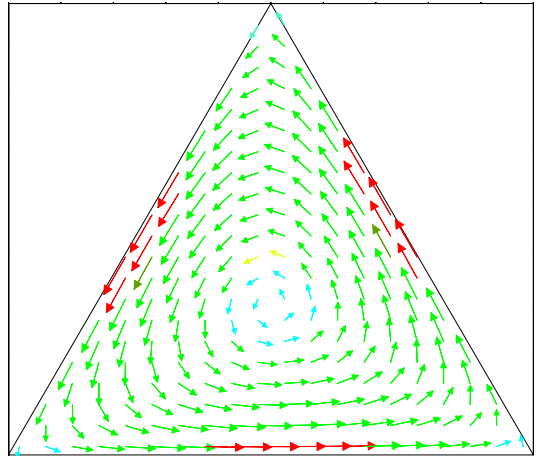

(a)

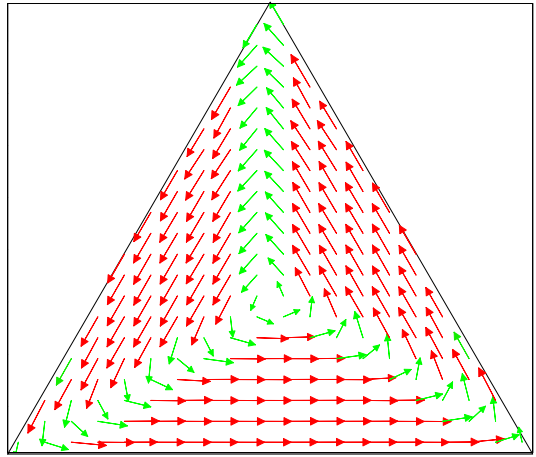

(b)

**Fig. 12.** Shear stress ($kN/cm^2$) vector distribution in the (a) elastic ($\theta/\theta_{el}=1$) and (b) elastic-plastic ($\theta/\theta_{el}=4.3$) regime in Example 4.2.

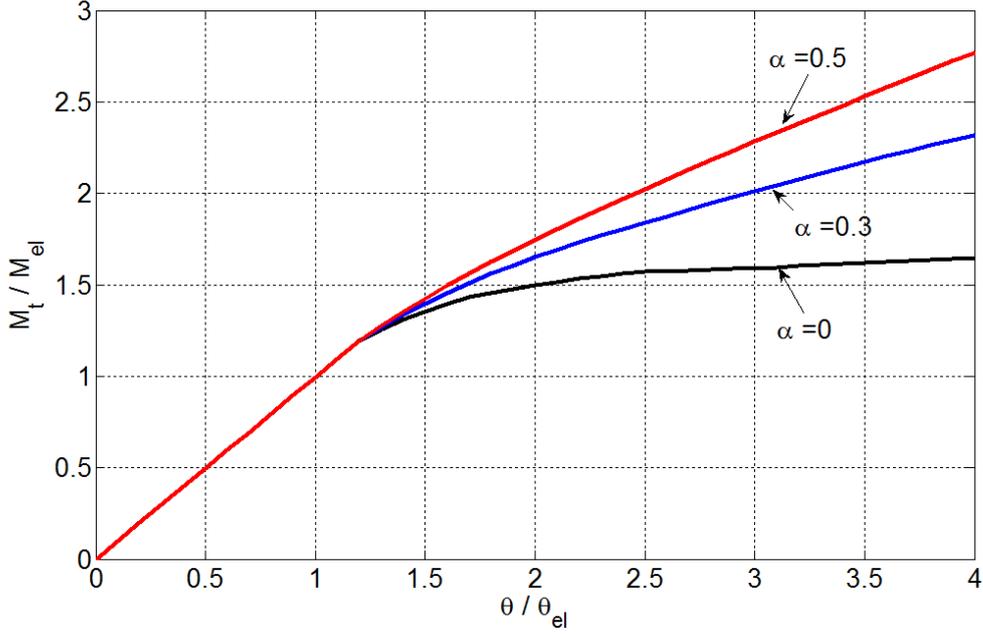

**Fig. 13.** Torsional moment-rotation curves for various values of the hardening parameter $a$ in Example 4.2.

*4.3 FGM Elastic-plastic bar with rectangular cross-section*

The final example is devoted to the study of an elastic-plastic bar with rectangular cross section $b \times h$ with $b = 5 cm$ and $h = 10 cm$ (see Fig. 4) made of an FGM with ceramic and metal constituents. The employed data of the ceramic phase are: $E_c = 5000 \text{ kN/cm}^2$, $\nu_c = 0.25$ and for the metal phase are: $E_m = 3000 \text{ kN/cm}^2$, $\nu_m = 0.25$ and $E_m^h = 500 \text{ kN/cm}^2$, $\sigma_Y = 5 \text{ kN/cm}^2$. The distribution of the overall values of $E$, $\sigma_Y$ and $E^h$ along the y axis are shown in Fig. 14. The results were obtained using $N = 300$ constant boundary elements and $M = 450$ nodal points (see Fig. 4). The maximum values of the elastic and fully plastic torsional moment for this cross-section made of metal are $M_{el} = 0.142 hb^2 \sigma_Y = 177.53 \text{ kNcm}$ and $M_{pl} = 0.0962(3h-b)b^2 \sigma_Y = 300.62 \text{ kNcm}$, respectively.

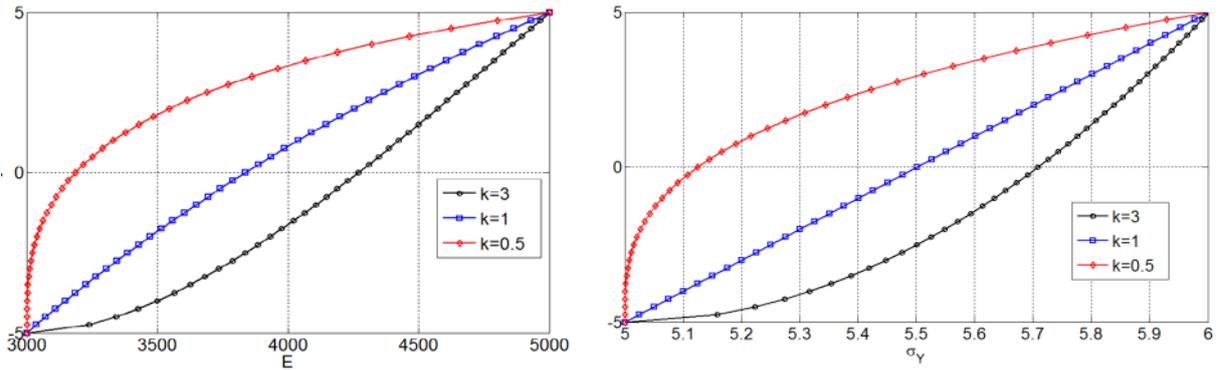

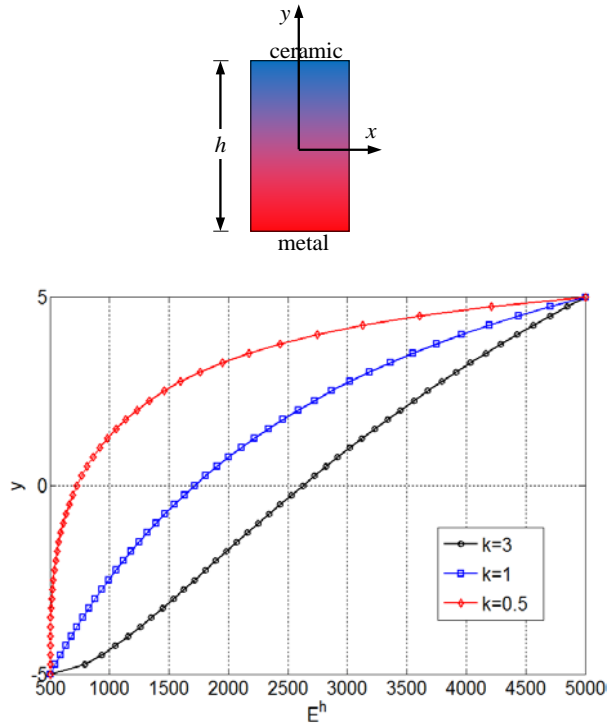

**Fig. 14.** Distribution of the overall values $E$, $\sigma_Y$ and $E^h$ along the y axis in Example 4.3.

Fig. 15 shows the torsional moment-rotation curve of the FGM material for various values of the power law exponent $k$. It is observed that the response of the FGM material approaches the response of the ceramic material when $k$ tends to zero, but for large values of $k$ approaches the response of the metal. Moreover, Fig. 16 - 20 present the progress of the plastic zones and shear stresses for various values of the power law exponent $k$ and the rotation $\theta$. What is interesting to notice is that the yielding always commences at points placed on the boundary independently of $k$. However, the plasticity spreads more rapidly as $k$ increases. Finally, the contours of the warping function in the elastic-plastic state ($\theta/\theta_{el} = 2.6$) are depicted in Fig. 21.

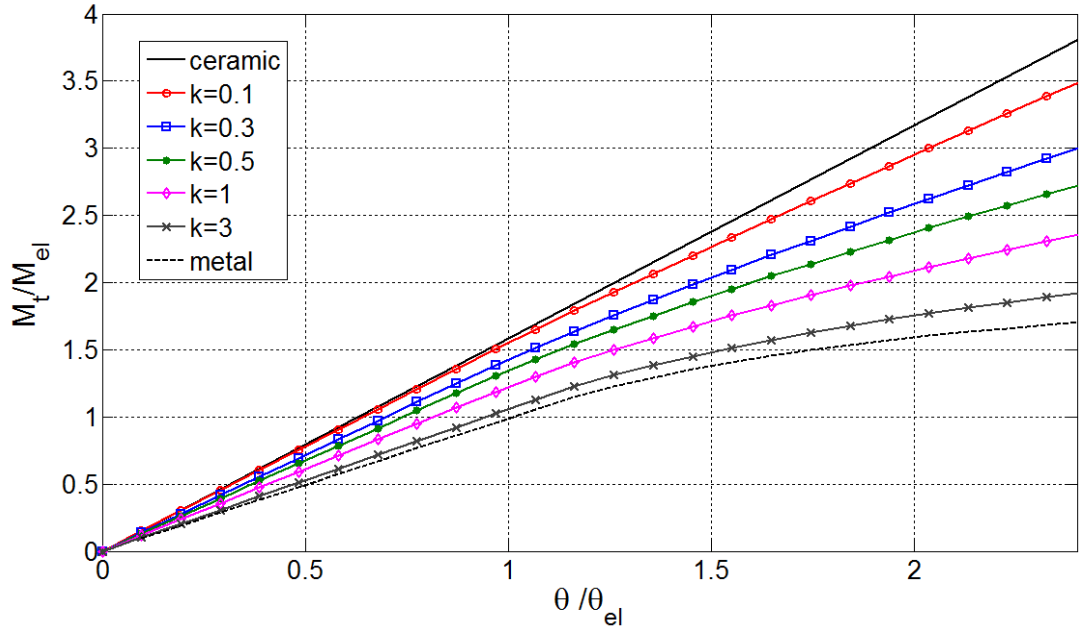

**Fig. 15.** Torsional moment-rotation curves of the FGM material for various values of the power law exponent $k$ in Example 4.3.

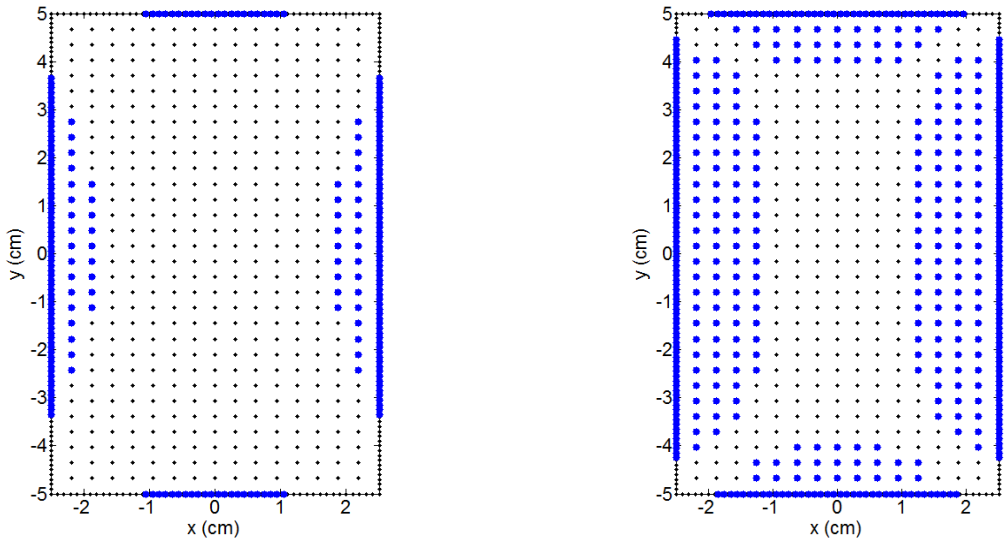

**Fig. 16.** Plastic regions for FGM material ($k = 0.1$) (a) $\theta/\theta_{el} = 1.06$ and (b) $\theta/\theta_{el} = 1.85$ in Example 4.3.

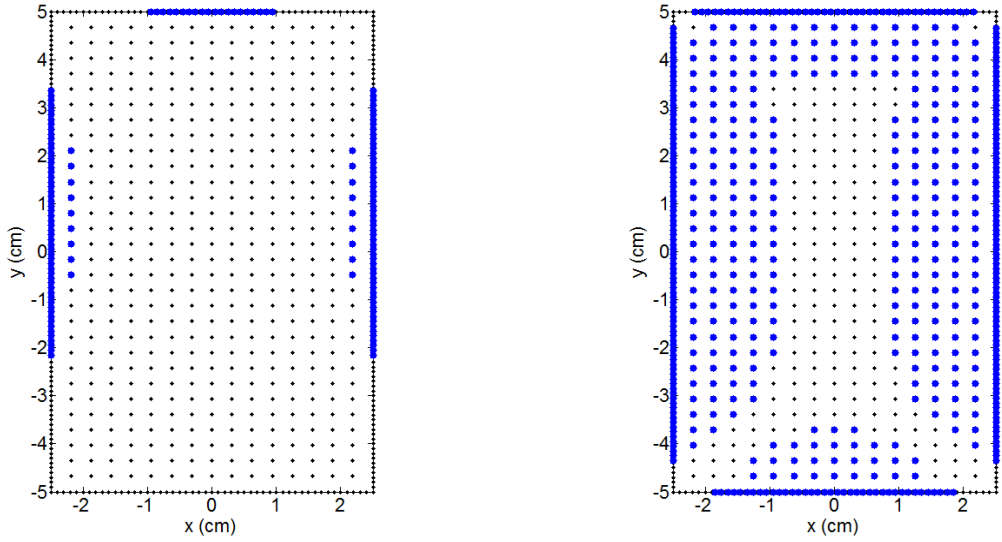

**Fig. 17.** Plastic regions for FGM material ($k=1$) (a) $\theta/\theta_{el}=1.06$ and (b) $\theta/\theta_{el}=2.60$ in Example 4.3.

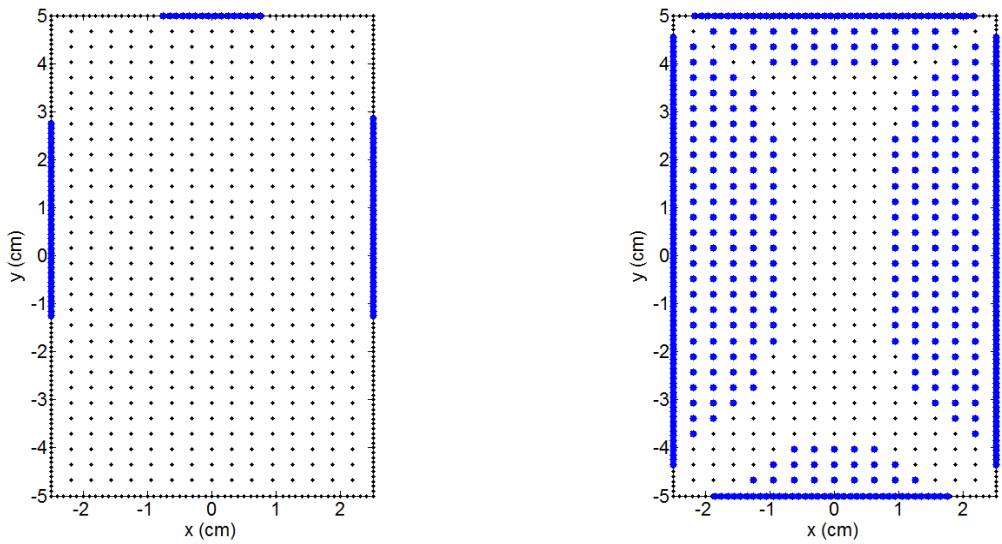

**Fig. 18.** Plastic regions for FGM material ($k=3$) (a) $\theta/\theta_{el}=1.06$ and (b) $\theta/\theta_{el}=2.60$ in Example 4.3.

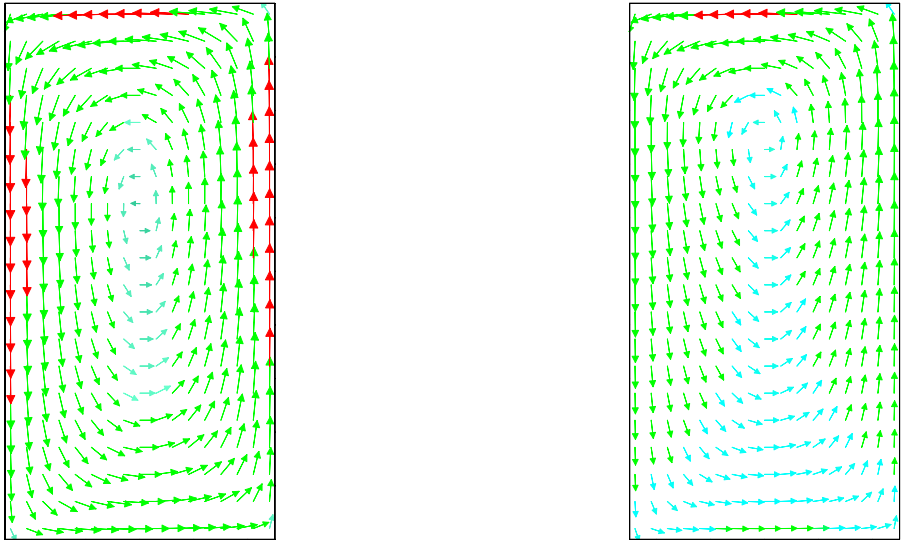

**Fig. 19.** Shear stress ($kN/cm^2$) vector distribution for FGM material ($k=1$) (a) $\theta/\theta_{el} = 1.06$ and (b) $\theta/\theta_{el} = 2.6$ in Example 4.3.

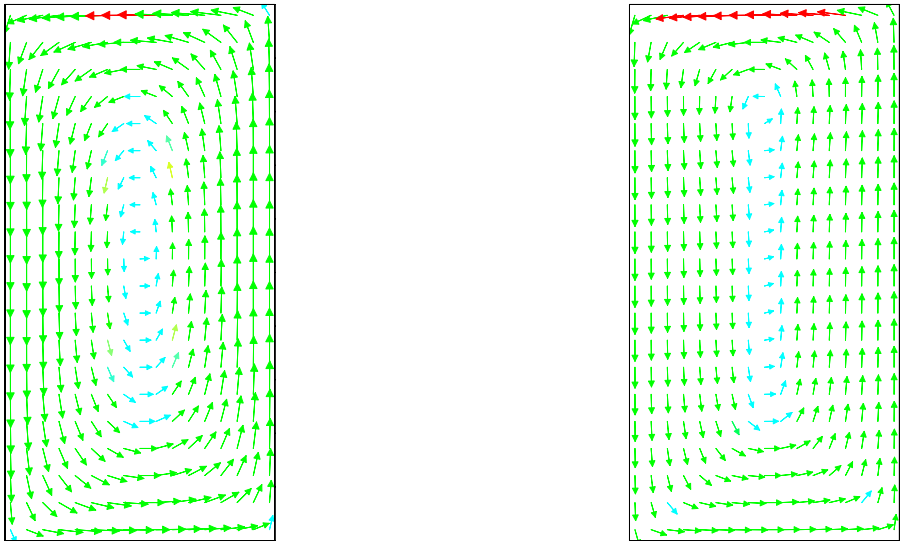

**Fig. 20.** Shear stress ($kN/cm^2$) vector distribution for FGM material ($k=10$) (a) $\theta/\theta_{el} = 1.06$ and (b) $\theta/\theta_{el} = 2.6$ in Example 4.3.

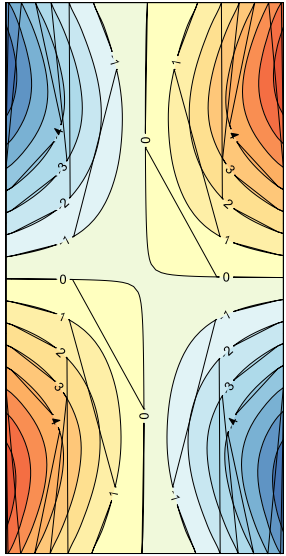

**Fig. 21.** Contours of the warping surface for FGM material ($k=1$) at $\theta/\theta_{el} = 2.6$ in Example 4.3.

## 5. Conclusions

In this paper a new integral equation solution to the elastic-plastic problem of functionally graded bars with arbitrary cross-section under torsional loading was presented. The nonlinear elastic-plastic problem was described mathematically by the deformation theory of plasticity. Several bars with various cross-sections and material types were analyzed, in order to validate the proposed model and exemplify its salient features. The main conclusions that can be drawn from this investigation are as follows:

1. The formulation is general in the sense that it can be applied to an arbitrary cross-section made of any type of elastoplastic material (elastic – perfectly plastic, strain hardening and nonlinear materials).
2. The proposed solution method is comprised of (i) an integral equation solution to the torsion problem of a non-homogeneous isotropic bar which has the following advantages:
   - The method provides a direct solution to the differential equation and overcomes the shortcoming of FEM solutions, which require resizing of the elements and re-computation of their stiffness during the iteration process,
   - The method is boundary only it has all the advantages of the BEM, i.e. the discretization and integration are performed only on the boundary,
   - The warping function and the stress resultants are computed at any point using the respective integral representation as mathematical formulae,

3. and (ii) a new straightforward nonlinear solution to the elastoplastic problem emanating from the deformation theory of plasticity with the following advantages:
   - The method is accurate and its implementation is simple. Moreover, it is alleviated from the drawback of existing iteration schemes (e.g. the projection method, the arc-length method, the Neuber's rule) which exhibit a certain degree of complexity.
   - The key characteristic of the proposed approach lies in the fact that a nonlinear system of algebraic equations is constructed and any numerical method can be employed for its solution.
4. For the elastoplastic torsion problem of homogeneous bars the deformation theory gives very good results as compared to those obtained by the flow theory and its main advantages are: simpler implementation and lower computational cost.
5. The torsional response of the metal-ceramic FGM approaches the response of the ceramic material when the power law exponent $k$ tends to zero, but for large values of $k$ approaches the response of the metal.
6. For the rectangular FGM bar is interesting to notice that the yielding always commences at points placed on the boundary independently of the power law exponent $k$. However, the plasticity spreads more rapidly as $k$ increases.
7. Overall, the ceramic constituent improves the torsional response of the FGM bar since it increases the plastic moment capacity of the cross-section and consequently its ultimate strength.

**References**

1. Mori, T. and Tanaka, K. (1973), "Average stress in matrix and average elastic energy of materials with misfitting inclusions," Acta Metall. 21, pp. 571–574.

2. Hill, R. (1965), "A self-consistent mechanics of composite materials," J. Mech. Phys. Solids 13, pp. 213–222.

3. Nakamura, T., Wang, T. and Sampath, S. (2000), "Determination of properties of graded materials by inverse analysis and instrumented indentation," Acta Mater. 48, pp. 4293–4306.

4. Tamura, I., Tomota, Y. and Ozawa, H. (1973), "Strength and ductility of Fe–Ni–C alloys composed of austenite and martensite with various strength," *Proceeding of the 3rd International Conference on Strength of Metals and Alloys*, Cambridge Institute of

Metals, Vol. 1, pp. 611–615.

5. Bocciarelli, M., Bolzon, G. and Maier, G. (2008), "A constitutive model of metal–ceramic functionally graded material behavior: formulation and parameter identification," Comput. Mater. Sci. 43, pp. 16–26.

6. Nadai, A. (1931), *Plasticity*, McGraw-Hill, New York.

7. Ponter, A.R.S. (1966), "On plastic torsion," Int. J. Mech. Sci. 8, pp. 227-235.

*8.* Mendelson, A. (1968), "Elastic–plastic torsion problem for strain-hardening materials," *NASA, TN D-4391*.

9. Hodge, P.G., Herakovich, C.T. and Stout, R.B. (1968) "On numerical comparison in elastic- plastic torsion," J. Appl. Mech.-T ASME 35, pp. 454-459.

10. Mendelson, A. (1975), "Solution of elastoplastic torsion problem by boundary integral method," *NASA, TN D-7872*.

11. Baba, S. and Kajita, T. (1982), "Plastic analysis of torsion of a prismatic beam," Int. J. Num. Meth. Engng. 18, pp. 927–944.

12. Zhen-Sheng, C. (1983), "A boundary element solution to elasto-plastic torsion of solids of revolution," Int. J. Num. Meth. Engng. 19, pp. 1193–1207.

13. Bilinghurst, A., Williams, J.R.L., Chen, G. and Trahair, N.S. (1992), "Inelastic uniform torsion of steel members," Comput. Struct. 42, pp. 887–894.

14. Wagner, W. and Gruttmann, F. (2001) "Finite element analysis of Saint-Venant torsion problem with exact integration of the elastic-plastic constitutive equations," Comp. Meth. Appl. Mech. Eng. 190, pp. 3831–3848.

15. Dwivedi, J.P., Shah, S.K., Upadhyay P.C. and Talukdar, N.K. Das (2002) "Springback analysis of thin rectangular bars of non-linear work-hardening materials under torsional loading," Int. J. of Mech. Sci. .44, pp. 1505-1519.

16. Sapountzakis, E.J. and Tsipiras, V.J. (2008), "Nonlinear inelastic uniform torsion of bars by BEM," Comput. Mech. 42, pp. 77–94.

17. Sapountzakis, E.J. and Tsipiras, V.J. (2011), "Inelastic nonuniform torsion of bars of doubly symmetric cross section by BEM," Comput. Struct. 89, pp. 2388–2401.

18. Kolodziej, J.A. and Gorzelanczyk, P. (2012), "Application of method of fundamental solutions for elasto-plastic torsion of prismatic rods," Eng. Anal. Bound. Elem. 36, pp.

81–86.

19. Zieniuk, E. and Boltuc, A. (2015), "The effective strategy for calculating strains and stresses in elasto-plastic torsion of a bar," J. Theor. Comput. Mech. 53, pp. 179-194.

20. Mukhtar, F.M. and Al-Gahtani, H.J. (2016), "Application of radial basis functions to the problem of elasto-plastic torsion of prismatic bars", Appl. Math. Model. 1, pp. 436–450.

21. Sapountzakis, E.J. and Tsipiras, V.J. (2009), "Nonlinear inelastic uniform torsion of composite bars by BEM," Comput. Struct. 87, pp. 151–166.

22. Sapountzakis, E.J. and Tsipiras, V.J. (2009), "Effect of axial restraint in composite bars under nonlinear inelastic uniform torsion by BEM," Eng. Struct. 31, pp. 1190–1203.

23. Bayat, Y. and Toussi, H.E. (2015) "Elastoplastic torsion of hollow FGM circular shaft," J. Comput. Appl. Res. Mech. Eng. 4, pp. 165-180.

24. Jones, R.M. (2009), *Deformation theory of plasticity*, Bull Ridge Publishing, Blacksburg, VA.

25. Katsikadelis, J.T. (2002), "The analog equation method. A boundary-only integral equation method for nonlinear static and dynamic problems in general bodies," Theor. Appl. 27, pp. 13–38.

26. Jahed, H., Sethuraman, R. and Dubey, R.N. (1997), "A variable material property approach for solving elastic–plastic problems," Int. J. Pres. Ves. Pip. 71, pp. 285–291.

27. Desikn, V. and Sethuraman, R. (2000), "Analysis of material nonlinear problems using pseudo-elastic finite element method," J. Press. Vess. Technol. 122, 457–461.

28. Dai, K.Y., Liu, G.R., Han, X. and Li, Y. (2006), "Inelastic analysis of 2D solids using a weak-form RPIM based on deformation theory," Comput. Methods Appl. Mech. Engrg. 195, pp. 4179–4193.

29. Katsikadelis, J.T. (2016), *The Boundary Element Method for Engineers and Scientists, 2nd Edition*, Academic Press, Elsevier, UK.

30. Katsikadelis, J.T. and Tsiatas, G.C. (2016), "Saint-Venant torsion of non-homogeneous anisotropic bars," J. Appl. Comput. Mech. 2, pp. 42-53.

31. Huang, H., Chen, B. and Han, Q. (2014), "Investigation on buckling behaviors of elastoplastic functionally graded cylindrical shells subjected to torsional loads," Compos. Struct. 118, pp. 234–240.

32. Katsikadelis J.T. (2008), "A generalized Ritz method for partial differential equations in domains of arbitrary geometry using global shape functions," Eng. Anal. Bound. Elem. 32(5), pp. 353–367.